\documentclass[12pt,a4paper]{article} 
\usepackage{amsfonts,amssymb} 
\usepackage[cp1251]{inputenc}

\usepackage[margin=25mm, a4paper]{geometry}  

\usepackage{tikz}  

\def\n{\noindent}  \def\?#1{}
\def\IZ{{\mathbb{Z}}}  \def\IR{{\mathbb{R}}}
 \def\cA{{\cal A}}   
\def\cT{{\cal T}}      \def\cS{{\cal S}}  
        \def\cN{{\cal N}}

   \def\phi{\varphi}  
\def\ep{\varepsilon}  \def\t{\tilde} \def\bdelete#1{}
    \def\v#1{\vec{#1}}  \def\f#1{\bar{#1}}

\def\mod1{\,({\rm mod\ } 1)\,}

\def\beq#1#2{\begin{equation} \label{#1} #2 \end{equation}}
\def\bea#1{\begin{eqnarray*} #1 \end{eqnarray*}} \def\a{\!\!\!&\!\!\!\!&}

\def\toas#1{\stackrel{#1}{\longrightarrow}}
\def\function#1{\left\{\!\!\!\begin{array}{ll} #1 \end{array} \right.}

\def\proof{\smallskip \noindent {\bf Proof. \ }}       
\def\blanksquare{\,\,\,$\sqcup\!\!\!\!\sqcap$}         
\def\qed{\hfill\blanksquare\linebreak\smallskip\par}   

\def\thname{Theorem}  \def\lmname{Lemma}    \def\prname{Proposition}  \def\qname{Question}
\def\dfname{Definition}  \def\crname{Corollary}  \def\rmname{Remark}
\def\exname{Example}  

\newtheorem{theorem}{\thname}[section]   
\newtheorem{lemma}{\lmname}[section]     
\newtheorem{proposition}[lemma]{\prname} 
\newtheorem{question}[lemma]{\qname} 

\newtheorem{dftn}{\dfname}[section]
\newenvironment{definition}{\begin{dftn}\rm}{\end{dftn}} 
\def\bdef#1{\begin{definition} #1 \end{definition}}
\newtheorem{rmrk}[lemma]{\rmname}
\newenvironment{remark}{\begin{rmrk}\rm}{\end{rmrk}}     

\makeatletter\def\fps@figure{htbp}\makeatother 

\begin{document}

\title{Inheritance of shadowing for dynamical semigroups}
\author{Michael Blank\thanks{
        Higher School of Modern Mathematics MIPT, 
        1 Klimentovskiy per., Moscow, Russia;}
        \thanks{National Research University ``Higher School of Economics'';
        e-mail: mlblank@gmail.com}
       }
\date{December 28, 2024} 
\maketitle

\begin{abstract} We extend the single-perturbation approach 
(developed in our earlier publications for the case of a single map) 
to the analysis of the shadowing property for semigroups of endomorphisms. 
Our approach allows to give a constructive representation for a true 
trajectory which shadows a given pseudo-trajectory.
One of the main motivations is the question of inheritance: 
does the presence of shadowing for all generators of a semigroup 
imply shadowing for the semigroup and vice versa. 
Somewhat surprisingly, the answer to these questions is generally negative. 
Moreover, the situation with shadowing turns out to be quite different 
in a semigroup and in a non-autonomous system, despite the fact that the 
latter can be represented as a single branch of the former. 
 \end{abstract}

{\small\n
}

\section{Introduction}
Chaotic dynamics is characterized by very fast (typically exponential) divergence of (initially) 
nearby trajectories. Therefore arbitrary small perturbations of chaotic dynamics 
(for example, round-off errors in numerical simulations) may change its behavior drastically. 
One of a very few theoretical justifications of numerical analysis of chaotic dynamics 
is based on the so called shadowing property. Namely, on the ability to trace trajectories 
of a weakly perturbed system by true trajectories of the original system. By means 
of the recently established single perturbation approach for shadowing (see \cite{Bl22,Bl22a}) 
the shadowing property was studied  for a broad family of endomorphisms. 

The ability to deal with dynamical systems with discontinuities is one of 
the major goals of this approach. Note that in our main results formulated 
in Theorems~\ref{t:generators}, \ref{t:main} we do not assume anything about 
the smoothness of the dynamical system. 
Examples of systems with discontinuities satisfying the single perturbation 
approach were discussed in detail in our earlier publications \cite{Bl22,Bl22a}.

The present paper is dedicated to the extension of this approach for non-autonomous 
discrete time systems and semigroups of endomorphisms (to which we refer as 
dynamical semigroups).

Recall that a non-autonomous discrete time system is a bi-infinite sequence of maps 
 $\v{g}:=\{g_n\}_{n\in\IZ}$ acting in the same space $X$. 
A (discrete time) semigroup $G$ of maps is a closed under 
composition collection of maps acting on the same space $X$. 

Informally, at time $t=n$ one applies the map $g_n$ in the first setting or 
an arbitrary element $g$ of the semigroup $G$ in the second setting. 

It is worth note, that in both settings one needs to modify the single perturbation 
approach to be applied at arbitrary time in distinction to the original definition for 
a single map, where only the perturbation at time $t=0$ was taken into account.

Given a point $x\in X$ one can define its trajectory under 
the action of the semigroup in two very different ways. First of them, 
to which we will refer as a global trajectory, is the union of all images of 
the point $x$ under the action of the maps from the semigroup.  
The main disadvantage of this definitions is that the points of the global 
trajectory are not ordered. To avoid this difficulty we consider another approach. 

One says that a semigroup is generated by a collection of maps $G:=\{g_i\}$ 
(called generators) if each element of the semigroup can be represented as 
a finite composition of the generators. From now on, we 
identify a semigroup to a given set of its generators, while the same 
semigroup can have very different sets of generators. In this terminology 
by a {\em trajectory} of the semigroup $G$ we mean a sequence 
$\v{x}:=\{x_k\}_{k\in\IZ} \subset X$ such that $Gx_k \ni x_{k+1}~~\forall k$. 
Here $Gx := \cup_i g_ix$.  

Similarly for a given metric $\rho$ in $X$ one defines a {\em pseudo trajectory} 
of the semigroup $G$ as a sequence $\v{y}:=\{y_k\}_{k\in\IZ} \subset X$ such that 
the sequence of values $\{\rho(Gx_k ,x_{k+1})\}_{k\in\IZ}$ (to which we refer as 
amplitudes of perturbations) is small in a certain sense (see discussion below). 
Here the functional $\rho(A,B)$ characterizes the distance between sets $A,B$.

Note that two different sets of generators of the same semigroup can lead to very 
different sets of (pseudo)trajectories. Therefore to justify our approach one needs to 
answer whether various shadowing properties depend on the choice of generators, 
which will be done in Theorem~\ref{t:generators}. 

From the point of view of shadowing, there is an important qualitative difference 
between ordinary dynamical systems and semigroups, consisting in the 
fact that in the latter case a pseudo-trajectory $\v{y}$ and a true trajectory $\v{x}$ 
approximating it can be generated by very different sequences of generators. 
On the one hand, this allows some flexibility in the construction of the approximation, 
but, on the other hand, the proof of the absence of shadowing becomes much more complicated. 

The problem of shadowing of pseudo-trajectories\footnote{Approximate trajectories of 
    a system under small perturbations, already considered by G.~Birkhoff~\cite[(1927)]{Bi} 
    for a completely different purpose.} 
by true trajectories of a hyperbolic dynamical system was first posed by 
D.V.~Anosov~\cite[(1967-70)]{An,An2} as a key step in analyzing the structural stability 
of diffeomorphisms. 
A similar but less intuitive approach, called ``specification'', was proposed at the same time 
by R. Bowen \cite[(1975)]{Bo}. Roughly speaking, both approaches check that the errors do not 
accumulate during the dynamics: each approximate trajectory can be uniformly traced 
by the true trajectory over an arbitrary long period of time. 
Naturally, this is of great importance in chaotic systems, where even a small error in the 
initial position can lead to a significant (exponentially fast) divergence of trajectories over time.

Further development (readers are referred to the main results, generalizations, and 
numerous references to monographs dedicated to this subject, such as \cite{Pi} and \cite{PS}, 
as well as the textbook \cite{KB}), has demonstrated a deep connection between 
the shadowing property and various ergodic characteristics of dynamical systems.
In particular, it was shown that for diffeomorphisms under some technical assumptions 
the shadowing property implies the uniform hyperbolicity. 
To some extent, this restricts the theory of uniform shadowing to an important but very special 
class of hyperbolic dynamical systems.

To overcome this restriction about 30 years ago a new concept of average shadowing was 
introduced in \cite{Bl88}, which allows to give a possibility to extend significantly the range of 
perturbations under consideration in the theory of shadowing, in particular to be able to deal 
with perturbations of Gaussian type, which are small only on average but not uniformly. 
The original idea under this concept was both to extend the range of perturbations and 
to be able to deal with non-hyperbolic systems. 
Although the first goal was largely achieved (see discussion of outstanding issues below), 
the second goal was not: only smooth hyperbolic systems were studied.
Nevertheless, this concept gave rise to a number of subsequent works in which various shadowing 
options were presented (see, for example, \cite{KKO,KO,LS,Sa,Sa2,WOC}) and the connections 
between them have been studied in detail.

The most challenging aspect of analyzing the shadowing property is dealing with the infinite 
number of possible perturbations to the original system. This makes the problem highly 
non-local, and it is desirable to find a way to reduce it to a single perturbation in order 
to gain better control over the approximation accuracy.

Following this idea, we recently developed in \cite{Bl22} a fundamentally new “gluing” construction, 
which consists of an effective approximation of a pair of consecutive segments of the true 
trajectories of an autonomous dynamic system. We showed in \cite{Bl22} that, given the gluing 
construction and some mild assumptions on autonomous dynamic systems, various shadowing-type 
properties can be proved. See exact definitions and details of the construction in Section~\ref{s:pre}. 


The paper is organized as follows. In Section~\ref{s:pre} we give general definitions 
related to the shadowing property and introduce the key tool of our analysis -- the gluing  
property. In Section~\ref{s:main} we formulate the main result -- Theorem~\ref{t:main}, 
which deduces various versions of shadowing from the gluing property. 
Sections~\ref{s:generators} and \ref{s:main} are dedicated to the proofs of the 
main results related to the shadowing properties of semigroups, 
while in Section~\ref{s:non-auto} we formulate and prove the corresponding statements 
for non-autonomous dynamical systems and discuss their similarities and distinctions 
to the semigroups. 
The remaining part of the paper is devoted to the analysis of the question of inheritance: 
does the presence of shadowing for all generators of a semigroup 
imply shadowing for the semigroup and vice versa. It is shown that this is absolutely 
not the case. 
In addition, we discuss whether it is possible to derive various shadowing properties 
if they are satisfied by ``similar'' (in particular, topologically conjugated) systems. 
A number of open questions are formulated.

\section{Setting and main result}\label{s:pre}
\bdef{A {\em semigroup} $G$ of mappings from a set $X$ into itself is a set of 
endomorphisms that is closed under finite superpositions.
A sub-collection of mappings $\{g_1,\dots,g_N\}\subset G$, satisfying the condition that each 
element  of the semigroup may be represented in the form %
$g = g_{i_n}\circ g_{i_{n-1}}\circ\dots\circ g_{i_1},  ~~ i_k \in\{1,\dots,N\},\quad k\in\{1,\dots,n\},$
is called the semigroup set of {\em generators}. }

To simplify notation in the sequel we identify a semigroup with its set of generators 
and write $Gx:=\cup_i g_ix$. It is worth mentioning that the same semigroup may be 
represented by very different sets of generators. We do not assume that the number of 
generators is finite, so in principle it is possible that $N=\infty$.

\bdef{A {\em trajectory} of the semigroup $G$ starting at a point $x\in X$ is a sequence 
of points $\v{x}:=\{\dots,x_{-2},x_{-1},x_0,x_1,x_2,\dots\}\subset X$, 
for which $x_0=x$ and $x_{i+1}\in Gx_i$ for all available indices $i$. 
The part of $\v{x}$ corresponding to non-negative indices 
is called the {\em forward (semi-)trajectory}, while the part corresponding to non-positive indices 
is called the {\em backward (semi-)trajectory}.}

In distinction to conventional dynamical systems (semigroups with a single generator) 
the initial point $x=x_0$ does not determine uniquely both the forward trajectory and 
the backward trajectory (in the case of a single non-bijective generator only the 
backward trajectory is not uniquely determined). In other words, for a given $x=x_0$ 
there might be arbitrary many admissible trajectories. Additionally the backward 
trajectory might be finite (if its ``last'' point has no preimages). 
In this case we are speaking only about available indices $i$ in the definition.

\begin{remark}
Unlike autonomous dynamic systems, where the introduction of a backward trajectory 
of an irreversible dynamical system looks somewhat peculiar and leads to multiple 
reverse branches, in the case of semigroups this demonstrates complete 
symmetry between the multiple branches of both forward and backward trajectories.
\end{remark}  

To deal with perturbations we assume that our semigroup acts on a complete  
metric space $(X,\rho)$ and extend this metric to the sets of subsets of $X$ as follows: 
$$ \rho(A,B) := \inf_{(a,B)\in A\times B}\rho(a,b) .$$ 
Note that this functional is not a metric, in distinction to e.g. the Hausdorff metric
$$\rho_H(A,B) := \min\{\sup_{a\in A} \inf_{b\in B}\rho(a,b), ~\sup_{b\in B} \inf_{a\in A}\rho(a,b)\} .$$ 

\bdef{A {\em pseudo-trajectory} of the semigroup $G$ is a sequence 
of points $\v{y}:=\{\dots,y_{-2},y_{-1},y_0,y_1,y_2,\dots\}\subset X$, for which the sequence of 
numbers $\{\rho(Gy_i,y_{i+1})\}$ for all available indices $i$ satisfies a certain ``smallness'' condition. 
We classify points corresponding to non-negative or non-positive indices as forward or backward 
pseudo-trajectories.}

Introduce the set of ``moments of perturbations'':
$$\cT(\v{y}):=\{t_i:~ \gamma_{t_i}:=\rho(Gy_{t_i}, y_{t_{i+1}})>0, ~i\in \IZ\}$$ 
ordered with respect to their values, i.e. $t_i<t_{i+1}~\forall i$. We refer to the 
amplitudes of perturbations $\gamma_{t_i}$ as {\em gaps} between consecutive segments 
of true trajectories. Throughout the paper we always assume that these gaps 
are uniformly bounded from above.

\bdef{For a given $\ep>0$ we say that a pseudo-trajectory $\v{y}$ is of 
\begin{itemize}
\item[(U)] {\em uniform} type, if $\rho(Gy_i,y_{i+1})\le\ep$ for all available indices $i$.
\item[(A)] {\em small on average} type, if 
    $\limsup\limits_{n\to\infty}\frac1{2n+1}\sum\limits_{i=-n}^n\rho(Gy_i,y_{i+1})\le\ep$. 
\item[(S)] {\em single perturbations} type, if the set $\cT(\v{y})$ consists of a single point. 
\end{itemize}}
If the backward pseudo-trajectory is finite, in the case (A) only positive indices $i$ are 
taken into account, which leads to one-sided sums 
$\frac1{n+1}\sum\limits_{i=0}^n\rho(Gy_i,y_{i+1})$ 
instead of two-sided ones.

A few words about history of pseudo-trajectories. 
Initially, G.~Birkhoff \cite{Bi} and much later D.V.~Anosov \cite{An} introduced the classical 
U-type pseudo-trajectories, which correspond to uniformly small errors during the dynamics. 
The A-type was proposed by M.~Blank \cite{Bl88,Bl22} in order to take 
care about perturbations small only on average. The difference between the versions in 
\cite{Bl88} and \cite{Bl22} is that the first one was unnecessarily strong due to purely technical 
arguments used in the proofs. In particular, for true Gaussian perturbations, the probability of 
this event was zero. The S-type was introduced to highlight the importance of analyzing the 
seemingly simple case of a single perturbation. 

Clearly $U,S\subset A$, but $S\setminus A\ne\emptyset$. 
To simplify notation we will speak about $\ep$-pseudo-trajectories, when the corresponding 
property is satisfied with the accuracy $\ep$. 

\bigskip

The idea of {\em shadowing} in the theory of dynamical systems boils down to the question,  
is it possible to approximate the pseudo-trajectories of a given dynamical system with true trajectories? 
Naturally, the answer depends on the type of approximation.

\bdef{We say that a true trajectory $\v{x}$ {\em shadows} a pseudo-trajectory $\v{y}$ 
with accuracy $\delta$ (notation $\delta$-shadows): 
\begin{itemize}
\item[(U)] {\em uniformly}, if $\rho(x_i,y_i)\le\delta$ for all available indices $i$.
\item[(A)] {\em on average}, if 
      $\limsup\limits_{n\to\infty} \frac1{2n+1}\sum\limits_{i=-n}^n\rho(x_i,y_i)\le\delta$. 
\item[(L)] {\em in the limit}, if $\rho(x_n,y_n)\toas{n\to\pm\infty}0$.
\end{itemize}}
If the backward pseudo-trajectory is finite, only positive indices $i$ are 
taken into account, which leads to the one-sided sum $\frac1{n+1}\sum\limits_{i=0}^n\rho(x_i,y_i)$.

Similarly to pseudo-trajectories the U-type shadowing was originally proposed by D.V.~Anosov \cite{An}, 
while the A-type was introduced\footnote{The reason is that pseudo-trajectories with large 
       perturbations cannot be uniformly shadowed.} by M.~Blank \cite{Bl22}.
Naturally, the types of pseudo-trajectories and the types of shadowing may be paired 
in an arbitrary way.

\bdef{We say that a semigroup $G$ acting on a metric space $(X,\rho)$ satisfies 
the {\em $(\alpha+\beta)$-shadowing property} (notation $G\in \cS(\alpha,\beta)$) 
with $\alpha\in\{U,A,S\},~ \beta\in\{U,A,L\}$ if $\forall\delta>0~\exists\ep>0$ such that 
each $\ep$-pseudo-trajectory of $\alpha$-type can be shadowed in the $\beta$ sense with the 
corresponding accuracy $\delta$.}

In this notation $\cS(U,U)$ corresponds to the classical situation of the uniform shadowing 
under uniformly small perturbations, while $\cS(A,A)$ corresponds to the average shadowing 
in the case of small on average perturbations.

One of the most interesting and open questions related to the shadowing problem is: 
under what conditions on a map does the presence of a specific type of shadowing for 
each pseudo-trajectory with a single perturbation imply a certain type of shadowing 
property for the system?
The reason for this question arises from the fact that the case of a single perturbation 
is simpler, and therefore, the idea of gaining insight into other types of perturbations 
from this knowledge is quite appealing. 
The answer is known, although only partially, in the case of U-shadowing for so-called 
positively expansive\footnote{Roughly speaking, expansivity means that if two forward 
    trajectories are uniformly close enough to each other, then they coincide. 
    In particular, this property is satisfied for expanding maps.}  
dynamical systems under the additional assumption that the single perturbation does 
not exceed $0<\ep\ll1$. 

To give the complete answer to this question we introduce the following property.

\bdef{We say that an $S$-type pseudo-trajectory $\v{y}$ with the perturbation at time $t=t_0$ 
is approximated by a true trajectory $\v{x}$ with accuracy rate $\phi:\IZ\to\IR_+$ {\em strongly} if 
\beq{e:glu}{
     \rho(x_k, y_k)\le\phi(k-t_0)\rho(Gy_{t_0-1},y_{t_0}) ~~\forall k\in\IZ}
and {\em weakly} if   
\beq{e:glu-w}{
     \rho(x_k, y_k)\le\phi(k-t_0) ~~\forall k\in\IZ .}  
}

This means that the true trajectory $\v{x}$ approximates both the backward and forward 
parts of the pseudo-trajectory $\v{y}$ with accuracy controlled by the rate function $\phi$, and 
in the strong version the accuracy additionally depends multiplicatively on the 
distance between the ``end-points'' of the glued segments of trajectories.

Without loss of generality, we assume that the functions $\phi(|k|)$ and $\phi(-|k|)$ 
are monotonic. Indeed, replacing a general $\phi$ by its monotone envelope 
$$\t\phi(k):=\function{\sup_{i\le k}\phi(i) &\mbox{if } k<0 \\  
                               \sup_{i\ge k}\phi(i) &\mbox{if } k\ge0} ,$$ 
we get the result. 


\bdef{We say that the semigroup $G$ satisfies the (strong/weak) 
{\em single perturbation approximation property} 
with the rate-function $\phi:\IZ\to\IR$ (notation $G\in \cA_{s/w}(\phi)$) if for each  
$S$-type pseudo-trajectory $\v{y}$ with $\cN(\v{y})=1$  there is a trajectory $\v{x}$ 
approximating $\v{y}$ in the strong/weak sense with accuracy $\phi$.}

\begin{remark} This property is a suitable reformulation for the present setting of the 
gluing property introduced earlier for conventional maps in \cite{Bl22,Bl22a}.
\end{remark}

\begin{remark}
In the case of a semigroup it is enough to check the that the inequalities (\ref{e:glu}) or (\ref{e:glu-w}) 
for the single perturbation at time $t_0=0$. 
\end{remark}
Indeed, for a given $\tau\in\IZ$ consider a trajectory $\v{x'}$ obtained from 
$\v{x}$ by the time shift by $\tau$, namely $x'_i:=x_{i+\tau}, ~\forall i$. 
Then assuming that $\v{y'}$ is approximated by $\v{x'}$ with accuracy $\phi(k), ~k\in\IZ$,  
we deduce the same property for $\v{x}, \v{y}$  with accuracy $\phi(k-\tau), ~k\in\IZ$. \qed

Our main results consist in the following statements.

\begin{theorem}\label{t:generators}
Let $G$ be a semigroup of maps from a complete metric space $(X,\rho)$ into itself. 
Then if for some choice of generators $G\in \cS(\alpha,\beta)$ 
with $\alpha\in\{U,A,S\},~ \beta\in\{U,A,L\}$, then the same property holds 
true for any other choice of generators. 
\end{theorem}

\begin{theorem}\label{t:main}
Let $G$ be a semigroup of maps from a complete metric space $(X,\rho)$ into itself. 
\begin{enumerate}
\item[(a)] If $G\in \cA_{s/w}(\phi)$ with $\phi(n)\toas{n\to\pm\infty}0$, then $G\in \cS(S,L)$.
\item[(b)]  If $G\in \cA_s(\phi)$ with $\Phi:=\sum_k\phi(k)<\infty$, then $G\in \cS(U,U) \cup \cS(A,A)$.
\end{enumerate}
\end{theorem}

In Section~\ref{s:non-auto} we will consider similar statements for non-autonomous systems.

\section{Proof of Theorem~\ref{t:generators}}\label{s:generators}
It can be shown that for a given semigroup, the minimal\footnote{A set of generators 
         is called minimal if no one of them can be represented as a superposition of 
         a finite combination of the others.} 
set of generators is uniquely determined up to a permutation. Nevertheless, the 
assumption of minimality is not necessary for the set of generators, and in the 
current setting of semigroups of mappings this assumption is very difficult to verify.
This explain importance of Theorem~\ref{t:generators}. 

\proof 
The main idea is that each new generator $h_j$ can be represented 
as a finite superposition of old ones $g_i$. Therefore the perturbations of the actions of 
the semigroup $G$ written in terms of the new generators $\{h_j\}$ inevitably occur 
much more rare in comparison to the representation in terms of $\{g_i\}$. 
On the other hand, the new presentation preserves the amplitudes of perturbations, 
which implies the result. \qed

\section{Proof of Theorem~\ref{t:main}}\label{s:main}

\proof The part (a) is a direct consequence of the single perturbation 
approximation property. 

The strategy of the proof of the part (b) follows basically the construction developed 
for similar statements in the case of conventional maps in \cite{Bl22,Bl22a}, 
with differences related to the more complex structure of the mapping semi-groups.
In particular, we prove that there is a constant $K=K(\phi)<\infty$, such that for each 
$\ep>0$ small enough there is a true trajectory uniform/average approximating with 
accuracy $\delta\le K\ep$ each $\ep$-pseudo-trajectory of U/A-type. 

Consider the set of moments of perturbations of a pseudo-trajectory $\v{y}:=\{y_i\}_{i\in\IZ}$ 
$$\cT(\v{y}):=\{t_i:~ \gamma_{t_i}:=\rho(Gy_{t_i}, y_{t_{i+1}})>0, ~i\in \IZ\} .$$ 
Between the moments of time $t_i$ there are no perturbations and hence $\v{y}$ can 
be divided into segments of true trajectories. Thanks to the $\cA_s(\phi)$ property each pair of 
consecutive segments of true trajectories considered as a pseudo-trajectory with a single 
perturbation can be approximated by a true trajectory with the controlled accuracy. 

Without loss of generality, we assume that perturbations occur at each time step 
and hence $t_i:=i~\forall i\in\IZ$.

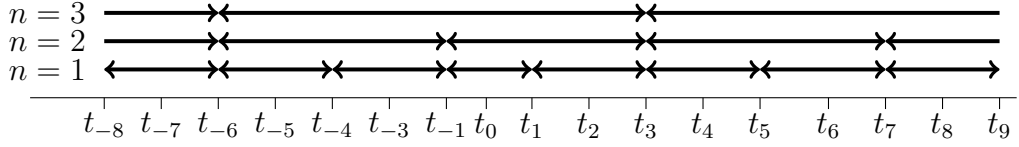
\begin{figure}\begin{center}
\begin{tikzpicture}[scale=0.75]
      \draw [-](-.8,0) to (16.5,0);
      \draw (.5,0) to (0.5,-.2); \node at (.5,-.5){$t_{-8}$};
      \draw (1.5,0) to (1.5,-.2); \node at (1.5,-.5){$t_{-7}$};
      \draw (2.5,0) to (2.5,-.2); \node at (2.5,-.5){$t_{-6}$};
      \draw (3.5,0) to (3.5,-.2); \node at (3.5,-.5){$t_{-5}$};
      \draw (4.5,0) to (4.5,-.2); \node at (4.5,-.5){$t_{-4}$};
      \draw (5.5,0) to (5.5,-.2); \node at (5.5,-.5){$t_{-3}$};
      \draw (6.5,0) to (6.5,-.2); \node at (6.5,-.5){$t_{-1}$};
      \draw (7.2,0) to (7.2,-.2); \node at (7.2,-.5){$t_{0}$};
      \draw (8.,0) to (8.0,-.2); \node at (8.0,-.5){$t_{1}$};
      \draw (9.,0) to (9.0,-.2); \node at (9.0,-.5){$t_{2}$};
      \draw (10.,0) to (10.0,-.2); \node at (10.0,-.5){$t_{3}$};
      \draw (11.,0) to (11.0,-.2); \node at (11.0,-.5){$t_{4}$};
      \draw (12.,0) to (12.0,-.2); \node at (12.0,-.5){$t_{5}$};
      \draw (13.2,0) to (13.2,-.2); \node at (13.2,-.5){$t_{6}$};
      \draw (14.2,0) to (14.2,-.2); \node at (14.2,-.5){$t_{7}$};
      \draw (15.2,0) to (15.2,-.2); \node at (15.2,-.5){$t_{8}$};
      \draw (16.2,0) to (16.2,-.2); \node at (16.2,-.5){$t_{9}$};
      \node at (-.5,.5){$n=1$}; \node at (-.5,1.0){$n=2$}; \node at (-.5,1.5){$n=3$};
      \draw [line width=1.5pt,,<->] (6.5,.5) to (8.0,.5); 
      \draw [line width=1.5pt,,<->] (8.0,.5) to (10.0,.5);
      \draw [line width=1.5pt,,<->] (10.0,.5) to (12.0,.5);
      \draw [line width=1.5pt,,<->] (12.0,.5) to (14.2,.5);
      \draw [line width=1.5pt,,<->] (14.2,.5) to (16.2,.5);
      \draw [line width=1.5pt,,<->] (4.5,.5) to (6.5,.5); 
      \draw [line width=1.5pt,,<->] (2.5,.5) to (4.5,.5); 
      \draw [line width=1.5pt,,<->] (.5,.5) to (2.5,.5); 
      \draw [line width=1.5pt,,<->] (6.5,1.0) to (10.0,1.0); 
      \draw [line width=1.5pt,,<->] (10.0,1.0) to (14.2,1.0); 
      \draw [line width=1.5pt,,<->] (2.5,1.0) to (6.5,1.0); 
      \draw [line width=1.5pt,,->] (.5,1.0) to (2.5,1.0); 
      \draw [line width=1.5pt,,<-] (14.2,1.0) to (16.2,1.0); 
      \draw [line width=1.5pt,,<->] (2.5,1.5) to (10.0,1.5); 
      \draw [line width=1.5pt,,->] (.5,1.5) to (2.5,1.5); 
      \draw [line width=1.5pt,,<-] (10.0,1.5) to (16.2,1.5); 
\end{tikzpicture}\end{center}
\caption{Order of the parallel gluing.}\label{f:sgluing} 
\end{figure}

We first approximate (see Fig.~\ref{f:sgluing}) pairs of segments around the moments 
of perturbations $t_i$ with even indices: $i_{\pm2k}$, obtaining longer 
segments of the true approximating trajectories. At each next step, we simultaneously  
proceed with consecutive pairs of already obtained segments. Consequently, at each step of the 
construction, we get a new pseudo-trajectory, consisting of half the number of segments of 
true trajectories with exponentially increasing lengths, but with larger gaps between them 
(compared to the original gaps). In the limit, we obtain an approximation of the entire pseudo-trajectory. 

To estimate approximation errors, we find the accuracy of the approximation of a pair of 
segments of true trajectories: $v_{-N^-}, v_{-N^-+1},\dots,v_{-1}$ and $v_0,v_1,\dots,v_{N^+}$. 
By the property $\cA_s(\phi)$ there exists a trajectory $\v{z}\subset X$ such that 

$$ \rho(v_k, z_k)\le\phi(k)\rho(Gv_{-1},v_0) ~~\forall k\in \{-N^-,\dots,N^+\}$$
Therefore 
$$ \sum_{k=-N^-}^{N^+}\rho(z_k, v_k) \le \rho(Gv_{-1},v_0) \sum_{k}\phi(k) 
    = \Phi\cdot \rho(Gv_{-1},v_0) .$$
There are several important points here: 
\begin{enumerate}
\item Each moment of perturbation $t_i$ is taken into account only once during the entire 
        approximation process.
\item The approximation accuracy depends only on the gap $\rho(Gv_{-1},v_0)$ between the 
        ``end-points'' of the segments of trajectories glued together.
\item The gaps between the end-points in the next step of the procedure may become larger 
         than the original gaps in $\v{y}$. 
\end{enumerate}

On the $n$-th step of the process of approximation of pairs of segments of true trajectories 
we get a pseudo-trajectory $\v{z}^{(n)}$ with a bi-infinite collection of gaps 
$\{\gamma_{t_i}^{(n)}\}_{i\in\IZ}$. 

Using the single perturbation approximation assumption, we obtain a recursive upper estimate for the gaps:
\beq{e:rec-est}{
  \gamma_{t_i}^{(n+1)} \le \gamma_{t_i}^{(n)} 
                                   + \phi_-^{(n)}\gamma_{t_{i-1}}^{(n)}  + \phi_+^{(n)} \gamma_{t_{i+1}}^{(n)} ,}
where $\phi_\pm^{(n)}=\phi(\pm2^n)$. Indeed, the lengths of the glued segments 
of trajectories on the $n$-th step of the procedure is equal to $2^n$ and 
$\phi_-^{(n)}\gamma_{t_{i-1}}^{(n)}$ is the upper estimate for the approximation error coming 
from the left, while $\phi_+^{(n)} \gamma_{t_{i+1}}^{(n)}$ is the upper estimate for the 
approximation error coming from the right. 

Rewriting (\ref{e:rec-est}) as follows:
$$ \gamma_{t_i}^{(n+1)} \le (\phi_-^{(n)} + \phi_+^{(n)}) \gamma_{t_i}^{(n)} 
             + \left(  (1 - \phi_-^{(n)} - \phi_+^{(n)}) \gamma_{t_i}^{(n)} 
                       + \phi_-^{(n)}\gamma_{t_{i-1}}^{(n)}  + \phi_+^{(n)} \gamma_{t_{i+1}}^{(n)}  \right) ,$$
we see that in the 1st term the factor $\phi_-^{(n)} + \phi_+^{(n)}$ vanishes with $n$, 
while the 2nd term corresponds to the averaging operator of type $v_i \to (1-a-b)v_i + av_{i-1} + bv_{i+1}$, 
the recursive application of which flattens a sequence $\{v_i\}$ to a constant. 

To make this argument precise, we need the following simple real analysis inequality, 
proved in \cite{Bl22,Bl22a}.

\begin{lemma}\label{l:inf-prod} For any sequence of nonnegative real numbers $\{b_k\}_{k\ge1}$ we have 
$$ \limsup_{n\to\infty}\prod_{k\ge1}^n(1+b_k) \le e^{\limsup\limits_{n\to\infty}\sum_{k=1}^n b_k} .$$
\end{lemma}

In the case of uniformly small perturbations $\gamma_{t_i}^{(0)}\le\ep~~\forall i\in\IZ$. 
Hence by Lemma~\ref{l:inf-prod} the upper bound of gaps 
$\gamma^{(n+1)}:=\sup_i \gamma_{t_i}^{(n+1)}$ may be estimated 
from above as follows
$$  \gamma^{(n+1)} \le  \ep \exp\left(\sum_k \phi(k) \right) = \ep e^\Phi .$$

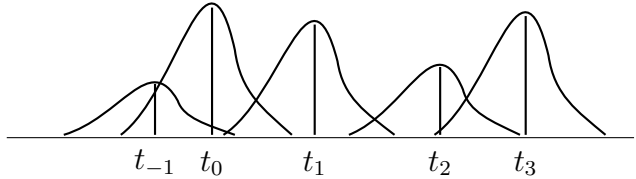
\begin{figure}\begin{center}
\def\G#1#2#3#4{\draw [thick] (1.6+#1,0*#3+#2) to (1.6+#1,1.5*#3+#2); \node at (1.6+#1,-.5+#2){#4};
      \draw [thick] (0+#1,0*#3+#2) .. controls (1+#1,.5*#3+#2) and (1.5+#1,2.5*#3+#2) .. (2+#1,1*#3+#2) 
                                 .. controls (2.1+#1,.5*#3+#2) and (2.5+#1,.3*#3+#2) .. (3+#1,0*#3+#2);}
\begin{tikzpicture}[scale=0.75]
      \draw [-](0,0) to (11,0); 
      \G{1}{.05}{.6}{$t_{-1}$}; \G{2}{.05}{1.5}{$t_{0}$}; \G{3.8}{.05}{1.3}{$t_{1}$};
      \G{6}{.05}{.8}{$t_{2}$}; \G{7.5}{.05}{1.4}{$t_{3}$};
\end{tikzpicture}\end{center}
\caption{Contributions to the upper bound of the gluing error.}\label{f:error} 
\end{figure}

Thus the gaps $\gamma_{t_i}^{(n)}$ are uniformly bounded from above by $\ep e^\Phi$. 
On the other hand, the contributions from each gap to the final approximation error 
are adding together (see Fig~\ref{f:error}). Using this, we get %
\beq{e:err-u}{  \rho(z_t^{(n)},y_t)  \le \ep e^\Phi  \sum_i \phi(t-t_i) \le \ep \Phi e^\Phi \quad\forall t\in\IZ . }

Similarly, for each $k>0$  one estimates the distance between pseudo-trajectories 
$\v{z}^{(n)}$ and $\v{z}^{(n+k)}$ as follows. Recall that the pseudo-trajectory $\v{z}^{(n)}$ 
consists of pieces of true trajectories of length $2^n$. Therefore 
$$ \rho(z^{(n)}_t,z^{(n+k)}_t)  \le \sum_{|j|\ge \tau(t,n)} \phi(j) \gamma_{j}^{(n)} ,$$
where $\tau(t,n)$ for a given $t$ grows to infinity as $2^n$. 

Since the function $\phi$ is summable, this implies that 
$$ \sum_{|j|\ge \tau(t,n)} \phi(j) \toas{n\to\infty} 0 .$$
Therefore, for any given $t$ the sequence $\{z^{(n)}_t\}_n$ 
is fundamental and converges as $n\to\infty$ to the limit $z_t$, where $\{z_t\}$ 
is the true trajectory of our system.

The estimate (\ref{e:err-u}) is uniform on $n$, so we may use it for $z_t$ as well, 
getting  
$$ \rho(z_t,y_t)  \le \ep\Phi e^\Phi \quad\forall t\in\IZ , $$
which proves that $G\in \cS(U,U)$.

\bigskip

Consider now the case of A-type perturbations. 
By the assumption that the perturbations are small on average, we get 
$$ \limsup_{k\to\infty}\frac1{2k+1}\sum_{i=-k}^k\gamma_{t_i}^{(0)} \le \ep .$$
Our aim is to show that $\exists C\ne C(\ep)$ such that 
$$ \limsup_{k\to\infty}\frac1{2k+1}\sum_{i=-k}^k\gamma_{t_i}^{(n)} \le C\ep \quad\forall n.$$

Without loss of generality, 
we assume that the function $\phi$ is even (i.e. $\phi(-k)=\phi(k)~~\forall k$). 
Indeed, replacing a general $\phi$ by its symmetrization $$\t\phi(k):=\max(\phi(-k), \phi(k))~~\forall k,$$ 
we get the result. 

Denote $R_k^{(n)}:= \sum_{i=-k}^k\gamma_{t_i}^{(n)}$. Then using (\ref{e:rec-est}) we get 

\bea{R_k^{(n+1)} \a= \sum_{i=-k}^k\gamma_{t_i}^{(n+1)} 
                             \le \sum_{i=-k}^k\gamma_{t_i}^{(n)} 
                                 + \sum_{i=-k}^k \phi_-^{(n)}\gamma_{t_{i-1}}^{(n)} 
                                 + \sum_{i=-k}^k  \phi_+^{(n)} \gamma_{t_{i+1}}^{(n)} \\ 
                            \a= R_k^{(n)} + \phi_-^{(n)} \left(\gamma_{t_{-k-1}}^{(n)} + R_k^{(n)} 
                                                                              - \gamma_{t_{k+1}}^{(n)} \right) \\
                           \a\qquad\quad+  \phi_+^{(n)} \left(-\gamma_{t_{-k-1}}^{(n)}  + R_k^{(n)} 
                                                            +  \gamma_{t_{k+1}}^{(n)}  \right) \\
                           \a= (1 + \phi_-^{(n)} + \phi_+^{(n)}) R_k^{(n)} 
                                  + (\phi_-^{(n)} - \phi_+^{(n)})(\gamma_{t_{-k-1}}^{(n)}  
                                                                                -  \gamma_{t_{k+1}}^{(n)}) \\
                           \a=  (1 + \phi_-^{(n)} + \phi_+^{(n)}) R_k^{(n)} 
                                          \qquad({\rm since}~\phi_-^{(n)} = \phi_+^{(n)}) \\
                           \a\le\dots \le R_k^{(0)} \prod_{i=0}^n  (1 + \phi_-^{(i)} + \phi_+^{(i)})  .} 

Contributions to the upper bound of the approximation error are coming from two different sources: 
the estimates of the gaps (changing during the steps of the parallel gluing construction) 
and the summation over error contributions from the approximation of pairs of consecutive segments 
of true trajectories (See Fig~\ref{f:error}).

To obtain an upper bound for the partial sum of gaps, setting $b_n:=\phi_-^{(n)} + \phi_+^{(n)}$, 
by Lemma~\ref{l:inf-prod} we get %
\beq{e:est-a} {R_k^{(n+1)} \le \prod_{i=1}^n(1+b_i) R_k^{(0)} \le e^{\sum_{i=0}^nb_i}R_k^{(0)} 
                        \le e^\Phi R_k^{(0)} .}%
Denoting $\gamma^{(n)}:=\sup_i \gamma_{t_i}^{(n)}$ and using (\ref{e:rec-est}), we get
$$ \gamma^{(n+1)}\le ( (1 + \phi_-^{(n)} + \phi_+^{(n)}))\gamma^{(n)} .$$
Applying Lemma~\ref{l:inf-prod} again, we obtain an upper estimate of gaps uniform on $n$:
\beq{e:bound-est}{\gamma^{(n)}:=\sup_i\gamma_{i}^{(n)}\le e^\Phi \gamma^{(0)}.}%

Now we are ready to estimate the average distance between the $n$-th approximating 
pseudo-trajectory $\v{z}^{(n)}$ and $\v{y}$. %
\bea{Q_k^{(n)} \a:= \frac1{2k+1} \sum_{t=-k}^k \rho(z_t^{(n)},y_t) \\
        \a\le  \frac1{2k+1} \sum_{t=-k}^k \sum_i \gamma_{-t+i}^{(n)} \phi(i) \\
          \a=  \sum_i \phi(i) \cdot \left( \frac1{2k+1} \sum_{t=-k}^k \gamma_{-t+i}^{(n)}\right)\\
          \a= \sum_i \phi(i) \cdot R_k^{(n)}(t), }
where $R_k^{(n)}(t):=\sum_{i=-k}^k\gamma_{-t+t_i}^{(n)}$.

Therefore, using (\ref{e:est-a}) we obtain an upper bound %
\beq{e:fin-a}{ \limsup_{k\to\infty} Q_k^{(n)}  \le \ep \Phi e^\Phi ,}%
which does not depend on the step number $n$. 

It remains to check that the sequence of approximating pseudo-trajectories $\v{z}^{(n)}$
converges to a limit $\v{z}$ being a true trajectory of our system. 
Exactly as in the case of uniformly small perturbations, we get
$$ \rho(z^{(n)}_t,z^{(n+k)}_t)  
     \le \sum_{|j|\ge \tau(t,n)} \phi(j) \gamma_{j}^{(n)} 
     \le e^\Phi \gamma^{(0)} \sum_{|j|\ge \tau(t,n)} \phi(j) \toas{n\to\infty}0 .$$
Thus for each $t\in\IZ$ the sequence $\{z_t^{(n)}\}_n$ is fundamental and converges to a limit $z_t$. 

Since the estimate (\ref{e:fin-a}) is uniform on $n$ we may use it for $\v{z}$ instead 
of $\v{z}^{(n)}$, getting  
$$ \limsup_{k\to\infty}\frac1{2k+1} \sum_{t=-k}^k \rho(z_t,y_t) \le \ep\Phi e^\Phi .$$
Theorem is proven. \qed

\section{Shadowing for non-autonomous systems}\label{s:non-auto}
A close look to the proof of Theorem~\ref{t:main} shows that a similar result may 
be obtained for a non-autonomous system considered as a single bi-infinite sequence 
of generators of a certain semigroup. 

To be precise let us give exact definitions.

\bdef{A {\em non-autonomous dynamical system} is defined by a time-dependent 
map $f:X\times \IZ\to X$, described by a bi-infinite sequence 
$\f{f}:=\{f_i\}_{i\in\IZ}$ of non-necessarily invertible maps $f_i:X \to X$ 
from a metric space $(X,\rho)$ into itself. In other words, $f(x,t):=f_tx,~t\in\IZ$.}

\bdef{A {\em trajectory} of the non-autonomous system $(\f{f},X)$ starting at a point $x\in X$ is 
a bi-infinite sequence of points 
$\f{x}:=\{\dots,x_{-2},x_{-1},x_0,x_1,x_2,\dots\}\subset X$, 
for which $x_0=x$ and $f_tx_t=x_{t+1}$ for all available indices $t$ (considered as time). }

\bdef{A {\em pseudo-trajectory} of the non-autonomous system $(\f{f},X)$ is a bi-infinite 
sequence of points $\v{y}:=\{\dots,y_{-2},y_{-1},y_0,y_1,y_2,\dots\}\subset X$, 
for which the sequence of distances $\{\rho(f_ty_t,y_{t+1})\}$ for all 
available indices $t$ satisfies a certain ``smallness'' condition. }

From the point of view of a semigroup with the set of generators $\{f_i\}$ this 
means that instead of considering all possible ``branches'' of the semigroup, 
we fix a single one, corresponding to a given bi-infinite sequence of generators 
representing $\v{f}$, and their perturbations in the case of a pseudo-trajectory. 

Further definitions of types of perturbations and shadowing follow the ones for 
the semigroup with the only difference that the bi-infinite sequence of generators 
representing the non-autonomous system $\f{f}$ is fixed. Despite appearances 
this implies a serious restriction to the definition of shadowing, namely 
instead of an arbitrary approximating trajectory of the semigroup one needs to 
consider only the ones generated by the given sequence of maps $\f{f}$. 

In other words, it is possible that a shadowing type result, which holds for the 
semigroup, may break down for certain non-autonomous systems 
inherited from this semigroup. 

Using exactly the same arguments as in the proof of Theorem~\ref{t:main} but 
applying them only for the given branch $\f{f}$ of the semigroup we get. 

\begin{theorem}\label{t:main-na}
The assumption $\f{f}\in \cA_\phi$ with $\sum_k\phi(k)<\infty$ implies 
$\f{f}\in \cS(U,U) \cup \cS(A,A)$.
\end{theorem}


\section{Inheritance of shadowing properties}\label{s:inheritance}
Numerous examples of the analysis of the shadowing properties for semigroups
with a single generator by means of the single-perturbation approach were discussed 
in detail in our earlier publications (see \cite{Bl22,Bl22a}). Therefore here we restrict 
ourselves to the peculiarities related to general dynamical semigroups and non-autonomous 
systems, and in particular their interconnections and inheritance properties. 

The table below demonstrates implications between pairs of types of perturbations 
and shadowing: 

\begin{center}
\begin{tabular}{cccccccccc}
     ~ & UU & UA & AU & AA & SU & SA & SL\\
   UU & =   & +   & -    & -    & +   & +  & ?  \\
   UA &  -? & =   & ?    & -?  & +   &  +  & ? \\
   AU &  +?& +?  & =   & +   & +   &  +  & ? \\
   AA & -   & -?   & -    & =   & -?  &  +  &-?\\
   SU & -? &   ?   & +?  & ?   & =   &  +  &+?\\
   SA & -  & +?   &   ?  & +?  & -   & =   &+?\\
   SL & -  &   ?   & -?   & -?   & -   &  ?  &= \\
\end{tabular}
\end{center}

The symbol ``+'' here indicates the true implication, while the symbol ``-'' means 
its absence. ``?'' means that the answer is not known, and ``$\pm$?'' denotes 
situations where the answer is only tentatively known. 
We see that serious additional work is required to fill in the gaps in this table.


\begin{proposition} \label{p:inv}
Let $f:(X,\rho) \to (X,\rho)$ be a bijection such that $\exists C<\infty$ for which %
\beq{e:u-dist}{ 
    \frac1C \rho(x,y) \le \rho(fx, fy)  \le C \rho(x,y)~~\forall x,y \in X }%
and let it satisfy a certain shadowing property. 
Then $f^{-1}:(X,\rho) \to (X,\rho)$ satisfies the same shadowing property. 
\end{proposition}
\proof 
Let $\v{x}:=\{x_k\}_{k\in\IZ} \subset X$ be a pseudo-trajectory of $f$ and let  
$\v{x'}:=\{x'_k\}_{k\in\IZ} \subset X$ be a shadowing it true trajectory. 
Consider a pair of sequences $\v{y}:=\{y_k:=x_{-k}\}_{k\in\IZ}$ 
and $\v{y'}:=\{y'_k:=x'_{-k}\}_{k\in\IZ}$. 

Then $\forall k\in\IZ$ 
$$ y'_{k+1} := x'_{-k-1} = f^{-1}x'_{-k} = f^{-1} y'_k ,$$
$$ \rho(y_{k+1}, f^{-1}y_k) = \rho(x_{-k-1}, f^{-1}x_{-k}) 
      = \rho(f^{-1} f x_{-k-1}, f^{-1}x_{-k}) \le \frac1C \rho(fx_{-k-1}, x_{-k}) .$$
Since $\v{x}$ is a pseudo-trajectory of $f$, the last inequality implies the claim. \qed

\subsection{Conjugation property and shadowability}\label{s:conj}

In what follows we will need several technical results related to the notion of 
topological conjugation. 

\bdef{A pair of maps $f:X\to X$ and $g:Y\to Y$ are said to be {\em topologically conjugated} 
if there exists a homeomorphism $h:X\to Y$ such that $h\circ f = g\circ h$.}

\begin{question} \label{q:conj-h} What are the necessary and sufficient conditions 
for topological conjugacy of homeomorphisms of the projectively extended real line?
\end{question}

It was expected that it would suffice to assume that they have the same number of fixed
points, but the following example (suggested by an anonymous referee) shows that this 
is not so. Indeed, the maps $x\to-x$ and $x\to2x$ have fixed points at $0$ and in $\infty$, 
but they are not topologically conjugate. At present the complete answer to this question 
is not known. However, for the convenience of readers, we present some interesting 
constructive examples of such conjugations.

Let $f(x):=x+1,~~g(x):=x+2$. Then the conjugating homeomorphism $h(x):=2x$. 
Thus the conjugation of two different shifts is even analytic. A very different situation 
can be found in hyperbolic maps, for example, the maps $f(x):=2x,~~g(x):=3x$
are conjugated by a non-Lipschitz homeomorphism $h(x):={\rm sign}(x) |x|^{\log_32}$. 

\begin{proposition} \label{p:conj}
Let $(G,X,\rho)$ and $(F,Y,\rho_Y)$ be two uniformly conjugated semigroups, 
namely $\forall g\in G, f\in F$ there is a homeomorphism $h_{g,f} :X\to Y$ such 
that $h_{g,f}\circ g = f \circ h_{g,f}$ and $\exists C<\infty$ for which %
\beq{e:u-dist0}{ 
    \frac1C \rho(x,x') \le \rho_Y(h_{g,hf}x, h_{g,hf}x')  \le C \rho(x,x')~~\forall x,x'\in Y,~~g\in G. }%
Then non-autonomous systems $(\v{g},X)$ and $(\v{f},Y)$ composed by 
the generators of the semigroups $F,G$ satisfy the same shadowing properties. 
If additionally %
\beq{e:u-dist1}{ 
    \frac1C \rho(x,x') \le \rho_Y(h_{g,f}x, h_{g',f'}x')  \le C \rho(x,x')~~\forall x,x'\in Y,
                                                                                                  ~~g,g'\in G,~~ f,f'\in F.} %
Then $(G,X,\rho)$ and $(F,Y,\rho_Y)$ satisfy the same shadowing properties. 
\end{proposition}
\proof
Let us start with the simplest situation with one-generator semigroups $F,G$ 
with conjugated generators $f:=h\circ g$ and ``conjugated'' pseudo-trajectories 
$\v{x}:=\{x_k\}_{k\in\IZ} \subset X$ and $\v{y}:=\{y_k:=h x_k\}_{k\in\IZ} \subset Y$. 
Assume that there exists the true trajectory $\v{x'}:=\{x'_k\}_{k\in\IZ} \subset X$ 
which shadows $\v{x}$ in a given sense. Then we claim that the sequence of points 
$\v{y'}:=\{y'_k:=h x'_k\}_{k\in\IZ} \subset Y$ shadows the pseudo-trajectory $\v{y}$ 
in the same sense. Indeed, 
$$ \rho_Y(y_k,y'_k) = \rho_Y(hx_k, hx'_k) \le C \rho(x_k, x'_k) ~~\forall k\in\IZ .$$
Since all shadowing conditions are based on distances between the corresponding 
points of the pseudo-trajectory and the true trajectory which shadows it, the last equality 
implies the result.

In the case of multiple generators the proof is a bit more tricky especially since the 
conjugating homeomorphisms differ for different pairs of generators. 
Let $\v{y}:=\{y_k:=\widehat{f_ky_{k-1}}\}_{k\in\IZ} \subset Y$ be a pseudo-trajectory 
of the semigroup $F$. Here the distances $\rho(f_ky_{k-1}, \widehat{f_ky_{k-1}})$ 
are small in the corresponding sense. Consider the  ``conjugated'' pseudo-trajectory 
$\v{x}:=\{x_k:=h_{g_k,f_k}^{-1}y_k\}_{k\in\IZ} \subset X$ of the semigroup $G$. 
By the assumption for $\v{x}$ there is a shadowing true trajectory 
$\v{x'}:=\{x'_k:=f_{k'}x'_{k-1} \}_{k\in\IZ} \subset X$. 
Note that the indices $k'$ and $k$ may differ from each other. 

Now the sequence $\v{y'}:=\{y'_k:=h_{g_{k'},f_{k'}}x'_k \}_{k\in\IZ} \subset Y$ 
is a trajectory of $F$, since 
$$ y'_k:=h_{g_{k'},f_{k'}}x'_k = h_{g_{k'},f_{k'}} f_{k'}x'_{k-1} = g_{k'} y'_{k-1} .$$
On the other hand, 
$$ \rho_Y(y_k,y'_k) = \rho_Y(h_{g_{k},f_{k}}x_k, h_{g_{k'},f_{k'}}x'_k) .$$

In the setting of non-autonomous systems $g_{k}=g_{k'}$ and $f_{k}=f_{k'}$. 
Therefore using (\ref{e:u-dist0}) we apply the same argument as in the case 
of single generators.

For general semigroups, the conjugating homeomorphisms  $h_{g_{k},f_{k}}$  
and $h_{g_{k'},f_{k'}}$ may be different, and to overcome this difficulty one 
additionally uses (\ref{e:u-dist1}) to obtain the result. \qed

\subsection{The simplest examples of shadowable and non-shadowable maps}

Consider a 4-parameter family $\Psi$ of 1-dimensional maps from $\IR$ into itself
$$ \psi_{a,b,c,d}(x):=\function{ax + c &\mbox{if } x\le0 \\  bx + d&\mbox{otherwise}.}$$
Here $a,b,c,d,x\in\IR$. 
The function $\psi$ will play the main building block in the constructions we will discuss in this Section. 

\begin{proposition} \label{p:p-lin}
The map $fx:=\psi_{a,b,c,c}(x)$ with $a,b>0$ belongs to the class $UU \cup AA$ iff $a,b<1$ or $a,b>1$.
\end{proposition}
\proof We start with the case $b=1$. Consider the forward part of the pseudo-trajectory $\v{x}$ 
of this map, defined by the relation $x_{n+1}:=f(x_n)+\ep$ with $x_0=1, n\in\IZ_+,~~\ep>0$. 
Then for any $y\in\IR$ the distance from the trajectory $f^n(y)$ to $\v{x}$ goes to $\infty$ 
as $n\to\infty$. This proves that $f(x):=\psi_{a,1,c,c}(x)$ is non-shadowable. 
In a similar way one, but considering the backward part of the pseudo-trajectory $\v{x}$, 
one studies the case $a=1$. 

Assume now that $a,b>1$. By Theorem~\ref{t:main} to prove the shadowing property 
it is enough to check that for any two semi-trajectories $f^{-n}u$ and $f^{n}v$ with $n\in\IZ_+$ 
the gluing property holds true, namely that there exists a point $w$ such that 
the true trajectory of this point $f^{n}w,~~n\in\IZ$ approximates both the mentioned 
semi-trajectories with exponentially fast decaying accuracy. Choosing $w:=v$, 
we get the result. Indeed, the forward part of the pseudo-trajectory $\v{v}$ coincides 
with $f^{n}v$, while the backward part converge to $f^{-n}u$ at rate $a^{-n}$.

In the case $0<a,b<1$ one uses the same construction, but with $w:=u$.

Consider different types of slopes $a<1<b$ or $a>1>b$ under the assumption . 
Let $a<1<b$ and $c=0$. Then the real line is divided into 3  $f$-invariant 
sets $(-\infty,0), \{0\}, (0,\infty)$. Moreover,
$$ f^n(x_0)\toas{n\to\infty}0, ~~ f^n(x_0)\toas{n\to-\infty}-\infty ~~\forall x_0<0, $$
$$ f^n(x_0)\toas{n\to-\infty}0, ~~ f^n(x_0)\toas{n\to\infty}\infty ~~\forall x_0>0. $$
Therefore an arbitrary small perturbation at point $x_0=0$ prevents the gluing property, 
which proves that the map with this set of parameters in not shadowable.

Similarly, one studies the case $a>1>b$ and $c=0$.

It remains to consider different types of slopes $a<1<b$ or $a>1>b$ but with $c\ne 0$. 
There are two very different possibilities: $c<0$ (which yields two fixed points) 
and $c>0$ (no fixed points). Despite these differences in both cases the gluing 
property cannot hold, but for different reasons. In the 1st case the backward and 
forward semi-trajectories staying at the fixed points cannot be glued together, 
while in the 2nd case any true trajectory diverges exponentially fast from 
at least one of backward and forward semi-trajectories starting from the points 
$\pm\ep$ (with arbitrary small $\ep\ne0$). \qed

\begin{remark}
A close look at the proof above reveals an apparent contradiction between the opposing 
shadowing properties for a number of parameter combinations for which the 
corresponding maps are topologically conjugated. 
This is explained by the additional assumption about the uniform conjugation in 
Proposition~\ref{p:conj}, which does not hold true (see examples just before 
Proposition~\ref{p:conj}).
\end{remark}

\begin{remark}
The previous results concern only maps of continuous type from the family $\Psi$, 
namely those with equal constant terms $c=d$ and positive slopes $a,b>$. If $c\ne d$, 
the corresponding maps become discontinuous at $x=0$, and different signs 
of the slopes lead to irreversibility. For various combinations of these parameters,
shadowing turns out to be an interesting question, and we leave them
as exercises for the reader.
\end{remark}

\subsection{Inheritance of shadowing from generators to semigroups and backwards}

At first glance, it seems that if all the generators of a semigroup $G$ satisfy a certain 
shadowing property, then the semigroup itself must also satisfy that property. 
Conversely, if the semigroup satisfies the property, then all its generators do as well. 
The goal of this section is to show that this is absolutely not the case. 
Moreover, we show that some nonautonomous dynamical systems composed of 
generators of a shadowable semigroup might be non-shadowable.
Therefore, for counterexamples we consider the most possible choice of semigroup 
generators, leaving the analysis of sufficient inheritance conditions for future research.

\begin{proposition} \label{p:p-gshad+Gnon}
Let $X:=\IR$, $g_1(x):=2x$, $g_2(x):=x/2$. Then $g_i\in \cS(U,U) \cup \cS(A,A)~\forall i$,
while the semigroup $G:=\{g_1,g_2\}\not\in \cS(U,U) \cup \cS(A,A)$.
\end{proposition}
\proof The first claim again follows from Proposition~\ref{p:p-lin}. 
To prove the second claim, consider a pair of one-sided semi-trajectories $g_2^{-n}u$ 
and $g_1^{n}v$ with $n\in\IZ_+$ of the semigroup $G$ 
and assume that $u-v$ is an irrational number. Then for any $y\in \IR$ each 
trajectory of the semigroup $G$, starting at the point $y$ at $t=0$, 
diverges exponentially fast from at least one of backward and forward semi-trajectories 
$g_2^{-n}u$ and $g_1^{n}v$ correspondingly. \qed

Finding an example of the opposite type when all $g_i$ are non-shadowable, 
while the semigroup $G:=\{g_i\}$ is shadowable turns out to be a difficult problem. 
In particular at present we do not have examples of this sort acting on the real line. 
A possible candidate was the pair of generators $g_1(x):=\function{x/2 &\mbox{if } x\le0 \\  2x&\mbox{otherwise}}$ and $g_2(x):=x+1$. Unfortunately at the moment we do not 
have a proof that $G:=\{g_1,g_2\}\in \cS(U,U) \cup \cS(A,A)$. The example below 
corresponds to the average shadowing only.

Let $X:=\{1,2,3\}$, $\rho(x,y):=\function{0 &\mbox{if } x=y \\  1&\mbox{otherwise}}$, 
$g(x):=(x+1)({\rm mod}3) + 1$.

\begin{proposition} \label{p:+-}
Let $g_1:=g$, $g_2:=g^{-1}$. Then $g_i\not\in \cS(A,A)~\forall i$, 
while the semigroup $G:=\{g_1,g_2\}\in \cS(A,A)$.
\end{proposition}
\proof The first claim was proved in \cite{Bl22a}. To prove that $G$ is shadowable
we apply Theorem~\ref{t:main}, namely we check the gluing property for an 
arbitrary pair one-sided backward and forward semi-trajectories $G^-u$ and 
$G^+v$ with $u\ne v\in X$ of the semigroup $G$. Now it is enough to find 
a true trajectory of the semigroup $G$, coinciding with the concatenation of 
the backward and forward semi-trajectories except for a few indices close to 0. 
This can be easily achieved having in mind that $g_1\circ g_2(x)=x~~\forall x$. 

To be precise, let us prove that $\forall n\ge3$ and $\forall u,v\in X$ there is $f\in G$, 
represented by a superposition of exactly $n$ generators of $G$, such that $v=fu$. 
Indeed, \\
-- if $u=v$, then $g_1g_2u=v$, $g_1g_1g_1u=v$, etc.\\
-- if $u+1=v ({\rm mod}3)$, then $g_1u=v$, $g_1g_1g_2u=v$, etc.\\
-- if $u+2=v ({\rm mod}3)$, then $g_2u=v$, $g_1g_1u=v$, $g_1g_2g_2u=v$, etc.\\
\qed

\begin{remark}\label{r:non-auto} 
The example considered in Proposition~\ref{p:+-} demonstrates 
also the striking difference between semigroups and non-autonomous systems 
composed by the same set of generators. The point is that the gluing in the 
proof above was achieved by exchanging certain generators in the construction 
of the approximating trajectory, which is not allowed in the setting of non-autonomous systems.
\end{remark}

\section{Acknowledgment}
The research is supported by the MSHE ``Priority 2030'' strategic academic leadership program.



\end{document}